\documentclass[11pt]{article}
\usepackage{amssymb}
\usepackage{amsmath}
\usepackage{setspace}

\begin{document}

\begin{center}
CHARACTERISTICS OF DISTANCE MATRICES, A SECOND LOOK
\end{center}

\vspace{0.5cm}

Bryan E. Cain

\vspace{0.5cm}

\noindent ABSTRACT

\vspace{0.5cm}

Here the definitions of nearest neighbor, robustness, concordance, and correlation, all of which feature in (Temple 2023) (henceforth abbreviated (T23)), are adjusted to make them completely mathematical while preserving their significance.  A characterization is given of the possible limits of a function of distance matrices as the data matrices from which they are derived acquire more and more columns while their number of rows and the distance defining norm (= coefficient, in the terminology of (T23)) are held fixed.  No data contribute to the discussion here, but many examples, with standard norms and data matrices having just a few rows and columns, play an important role.  Indeed, small data matrices are displayed showing that robustness, defined either of the two ways, can be zero.  

\vspace{0.5cm}

\noindent INTRODUCTION

\vspace{0.5cm}

(T23) compares different mathematical approaches to classification by applying a rich mix of analytical, statistical, geometrical, and probabilistic insights and methods to data sets--some empirical and some computer generated.   Although the following commentary on (T23) does not engage with any data, it does adduce many examples. These examples may initially seem too simple and artificial to interest the classification community, but these uncomplicated matrices keep the discussion concrete and vivid and chasten any who might over-believe in the role of mathematics in classification.  

Distance matrices are still at the heart of the discussion, as are again their robustness (defined two ways), concordance, and correlation. In choosing their formal definitions, we part company with other reasonable possible definitions.   Since we do not comment on all of (T23) there remains much for others to discuss.

We solicit our readers' pardons for having adopted some verbiage and notations differing from that in (T23).

\vspace{0.5cm}

Let $M(n \times k)$ denote the set of all $n \times k$ real matrices, and let $M(n)$ denote the set of all real matrices having $n$ rows and a finite number of columns.  Let $0(n \times k)$ denote the member of $M(n \times k)$ whose every entry is zero.  In order to define the distance between two members of $M(1)$, which have the same number of columns, we rely on norms. (There are, however, other ways to define distance.  Our word ``norm'' corresponds to (T23)'s word ``coefficient''.) We call a function $N$ defined on $M(1)$ a norm iff (= if and only if) for every $v$ and $w$ in $M(1)$ and every real number $c$:  
\begin{itemize}
\item[(i)] $N(v)$ is zero if $v = 0(1 \times k)$ for some $k > 0$, and $N(v)$ is positive otherwise;  
\item[(ii)] $N(cv) = |c|N(v)$; and
\item[(iii)] $N(v + w) \leq N(v) + N(w)$ whenever $v$ and $w$ have the same number of columns so that $v + w$ is defined.  (This is the triangle inequality: The length $N(v + w)$ of the side $v + w$ is dominated by the sum of the lengths of the other two sides $v$ and $w$.)
\end{itemize}

Though not part of the definition of norm, here are two additional properties the norms we consider will satisfy.  
\begin{itemize}
\item[(iv)] $N(v) = N(v')$, if $v \neq 0(1 \times k)$ for any $k > 0$ and $v'$ is $v$ with all its zero entries removed.
\item[(v)] $N([1]) = 1$.
\end{itemize}

When $N$ is a norm ``the open $k$-dimensional $N$-ball of radius $r > 0$ centered at $z \in$ $M(1 \times k)$'' is $\{w \in M(1 \times k)\ :\ N(w - z) < r \}$.
Below we shall introduce the $p$-norms, a class of norms $N$ satisfying (i)--(v) with $N(v)$ an explicit formula in the entries of $v$, and we turn to them whenever $v$ is specified numerically and we need to compute its norm.

To an $X$ in $M(n \times k)$ with $i$-th row denoted $x(i)$, a norm $N$ associates $D(N,X)$, an $n \times n$ so-called ``distance matrix".  The entry in row $i$ and column $j$ of $D(N,X)$ is defined to be $N(x(j) - x(i))$, the distance between $x(j)$ and $x(i)$ determined by $N$.  Then the entries of $D(N,X)$ are all nonnegative, are all zero on its diagonal, and it itself is symmetric:  $N(x(j) - x(i)) = N((-1)(x(i) - x(j))) = |-1|N(x(i) - x(j))$.  

By ``constant column'' we mean any member of $M(n \times 1)$ all of whose entries are the same, e.g. $0(n \times 1)$.  If $Y$ is obtained from $X$ by adding to its columns constant columns and adjoining to it $j - k > 0$ constant columns so that $Y$ is in $M(n \times j)$ then $D(N,Y) = D(N,X)$ because of (iv).
 
(This does not exhaust the set of $Y$ in $M(n)$ such that $D(N,Y) = D(N,X)$.  For example let $N$ be a norm on $M(1)$ and let $w \to w'$ be an origin fixing $N$-isometry of $M(1 \times k)$, that is $0(1 \times k)' = 0(1 \times j)$ for some integer $j > 0$ and $N(w') = N(w)$ for all $w \in M(1 \times k)$. By (Mazur and Ulam, 1932) this isometry is linear. Then $D(N,Y) = D(N,X)$ if $Y \in M(n \times j)$ and its $i$-th row  $y(i) = x(i)'$, because $N(y(j) - y(i)) = N(((x(j) - x(i))') = N(x(j) - x(i))$.  If another norm $M$ is involved in the discussion one may ask that $D(M,Y) = D(M,X)$ also.)

When constructing examples of distance matrices, it may be helpful to keep in mind that if $Y$ is $X$ with its $i$-th row removed then $D(N,Y)$ is $D(N,X)$ with its $i$-th row and column removed.  Furthermore, if $P$ in $M(n \times n)$ is a permutation matrix, i.e. is obtained from the identity matrix of order n by shuffling its rows, then $PX$ is the same shuffle of the rows of $X$ and $D(N,PX) = PD(N,X)P^{T}$ where $T$ stands for ``transpose''.  In this way entries in $D(N,X)$ can be repositioned by starting with $PX$ instead of $X$.

Let $h(D_{1}, ..., D_{m})$ be a real valued function of $m$ distance-matrix variables $D_{i}$.  For $ i = 1, \dots, m$ let $N_{i}$ be a norm.  Consider the real-valued function of $m$ variables $X_{i}$ in $M(n)$ given by

$$H(X_{1}, \dots, X_{m}) = h(D(N_{1},X_{1}), \dots, D(N_{m},X_{m})).$$

The apparent greater generality obtained by allowing the $X_{i}$ to have differing numbers of columns is illusory:  Just adjoin to each some constant columns e.g. copies of $0(n \times 1)$ to endow them all with the same number of columns.  So we may and will assume without loss that each $X_{i}$ is in
$M(n \times k)$ for some $k$.

What about the limit of $H$ as $k$ approaches $\infty$?  If $H$ is a constant function, then the limit is that constant.  If $H$ is not constant and it is conjectured $a$ is its limit then there is a $b = H(W_{1}, \dots, W_{m})$ with all $W_{i}$ in $M(n \times k)$ and $b \neq a$.  As we explained above, adjoining $j$ zero columns to each $W_{i}$ to get $V_{i}$ in $M(n \times (j + k))$ will leave the value of $H$ unaltered:  $b = H(V_{1}, \dots, V_{m})$.  Increasing the number of columns this way does not alter $b$, the value of $H$.  Hence, $H$ is not approaching $a$ as the number of columns increases.  Thus:

\vspace{0.5cm}

\noindent THEOREM:  $H$ has a limit as the number of columns increases iff $H$ is constant.

\vspace{0.5cm}

\noindent REMARKS

\begin{itemize}
\item[(a)] We have held $n$ constant while increasing the number of columns.  This seems consistent with (T23) and indeed, though (T23) experiments with different fixed values of $n$, there is no discussion of how $n$ and $k$ might move together, perhaps in some kind of coordination, as $k$ increases.  So we conclude that (T23) had in mind an arbitrary positive fixed $n$ as $k \to \infty$ that is, $k$ approaches $\infty$.

\item[(b)] This theorem challenges the remarks in Table 3 of (T23) which offer conclusions about limits as $k \to \infty$.
\end{itemize}

To escape this disappointing theorem, we try taking limits of $H$ along sequences.  Let $M(n)^{m=}$ denote the set of $m$-tuples of members of $M(n)$, every entry in each $m$-tuple having the same number of columns. A sequence in $M(n)^{m=}$ will be called ``strictly increasing'' iff the common number of columns in a member of the sequence is strictly less than the common number in the next sequence member.    Evaluate $H$ at each member of the strictly increasing sequence and try to take the limit of the resulting sequence of values of $H$.  A modest gain results:

\vspace{0.5cm}

\noindent THEOREM:  The collection of limits of $H$ along strictly increasing sequences is the topological closure of the range of $H$.  (That closure is the range together with the limits of all sequences of members of the range.)

\vspace{0.5cm}

\noindent PROOF:  If $a$ is a limit of the values $H$ takes along a strictly increasing sequence in $M(n)^{m=}$ then that sequence of values is a sequence of members of the range of $H$ so $a$ is also a limit point of the range of $H$.

On the other hand, if $a$ is a limit point of the range of $H$, there is a sequence in $M(n)^{m}$ and the sequence of values $H$ takes along this sequence converges to $a$.  By adjoining zero columns as required to the $m$ matrices in each member of the sequence we can produce a new $m$-tuple each of whose members has the same number of columns and at which $H$ takes the same value it took at the original $m$-tuple. Adjoining additional zero columns if necessary we can produce a strictly increasing sequence along which $H$ takes the same sequence of values it took along the original sequence.  Thus $a$ is a limit of $H$ along a strictly increasing sequence in $M(n)^{m=}$. ///

\vspace{0.5cm}

Here are some examples of $n \times k$ real matrices $X = [x_{ij}]$  with $i$-th row $x(i)$ and distance matrices $D(N,X)$ the norm $N$ associates with some of them.

Ex0:  If $n = 1, D(N,X) = [0]$, which is $1 \times 1$.

Ex1:  If $n = 2$ and $c = N(x(2) - x(1))$ then $D(N,X) = \begin{pmatrix} 0 & c\\c & 0 \end{pmatrix}$.

Ex2:  If $k = 1, D(N,X) = N([1])[|x_{j1} - x_{i1}|] =  [|x_{j1} - x_{i1}|]$.

Ex3:  If $w \in M(1 \times k)$ and $a = [a_{1} \  \dots \  a_{n}]^{T}$ and if $X = aw$, then $D(N,X) = N(w)[|a_{j} - a_{i}|]$.  Ex2 is a special case of this.

Ex4:  $$X = \begin{bmatrix} 0 & 0 \\ 2 &  0 \\  -1 & \sqrt{3} \\ -1 & -\sqrt{3} \end{bmatrix}$$

\noindent The rows are the origin and the vertices on the circle of radius $2$ of an equilateral triangle with one vertex at $(2, 0)$.  Let $N((x,y)) = \sqrt{x^{2}+ y^{2}}$ and $s = \sqrt{3}$ then

$$D(N,X) =  \begin{bmatrix} 0 & 2 & 2 & 2 \\ 2 & 0 & 2s & 2s \\ 2 & 2s & 0 & 2s \\ 2 & 2s & 2s & 0 \end{bmatrix}$$

Ex5:

$$X = \begin{bmatrix} 0 & 0 \\ 1 &  1 \\  -1 & 1 \\ -1 & -1 \\ -1 & 1 \end{bmatrix}$$

\noindent The rows are the origin and the vertices of an origin-centered square of edge 2 and sides parallel to the coordinate axes.  Let $N$ be as in Ex4 and set $r = \sqrt{2}$.  Then

$$D(N,X) =  \begin{bmatrix} 0 & r & r & r & r \\ r &  0 & 2 & 2r & 2 \\ r & 2 & 0 & 2 & 2r \\ r & 2r & 2 & 0 & 2 \\ r & 2 & 2r & 2 & 0 \end{bmatrix}.$$

\vspace{0.5cm}

Here are some widely used norms which satisfy (iv) and  (v).  Let $w = [w_{1} \  \dots \  w_{k}] \in M(1)$ and let $p \in J = [1, \infty]$.  Set

$$N_{\infty}(w) = \max\{|w_{i}| : i = 1, \dots, k\}$$

\noindent and when $p$ is finite, $N_{p}(w)$ is the nonnegative $p$-th root of

$$|w_{1}|^{p} + \cdots + |w_{k}|^{p}.$$

\noindent If $q \in J$ and $q < p$ then $N_{p} \leq N_{q}$.

These are the well-known $p$-norms.  In (T23) only the most important are considered: $p = 1, 2, \infty$ and these are referred to there respectively as Manhattan, Euclidean, Hausdorff and abbreviated M, E, H. (Comparing this Hausdorff distance with Hausdorff's distance between sets does not repay the effort.) 
We shall continue in sequel to let $N$ denote an arbitrary norm satisfying (iv) and (v).

\vspace{0.5cm}

\noindent NEAREST NEIGHBORS

\vspace{0.5cm}

As above $X$ is in $M(n \times k)$ and its $i$-th row is $x(i)$ and $N$ is a norm on $M(1)$ satisfying (iv) and (v). We now define ``$u$ is an $N$-nearest neighbor of $v$'', but only when both $u$ and $v$ are rows of an $X$ in $M(n \times k)$.

When $n = 1$ then $D(N,X) = [0]$ and we say: $x(1) = X$ has no nearest neighbor.  When $n > 1$ and $i$ and $j$ are unequal members of $\{1, ..., n\}$ we say:  row $x(j)$ is an ``$N$-nearest neighbor'' of row $x(i)$ iff $N(x(j) - x(i))$ is the smallest off-diagonal entry in row $i$ of $D(N,X)$. There can be ties, so
this may happen for more than one value of $j \neq i$.  (Note: If $x(j) = x(i)$ for differing $i$ and $j$ (i.e. the rows of $X$ form a multiset not a set.) then row $i$ of $D(N,X)$ will be 0 in column $j$ and 0 is the smallest possible off-diagonal entry.  Thus this definition calls $x(j)$ an $N$-nearest neighbor of $x(i)$ (and vice versa too by symmetry).  Some may not like this and would have preferred ``smallest positive off-diagonal entry...''.  That would lead to no row of $X$ having a nearest neighbor, if $x(1) = \cdots =
x(n)$.)

\vspace{0.5cm}

\noindent EXAMPLES

\vspace{0.5cm}

\noindent (1) In Ex1, the two rows are $N$-nearest neighbors of each other.  

\noindent (2) In Ex3, let $n = 3, a = [1 \  3 \  4], c = N(w) > 0.$   Then

$$D(N,X) = c\begin{pmatrix} 0&2&3\\2&0&1\\3&1&0 \end{pmatrix}.$$
 
\noindent Thus row 2 is the $N$-nearest neighbor of row 1 but not vice versa, and rows 2 and 3 are $N$-nearest neighbors of each other.  This does not surprise since all the rows are on the line through $w$ and the origin with their spacings given by their coefficients.

\noindent(3) In Ex4 and Ex5, the $N$-nearest neighbors can be easily read off from the displayed distance matrices of the geometry of the set of points since there the norm is the Euclidean $2$-norm.

\vspace{0.5cm}
           
Set $NEAR(N,X,i) = \{j \neq i : x(j)$ is a $N$-nearest neighbor of $x(i)\}$.

\vspace{0.5cm}
 
When $n > 1$, this set contains at least one integer and at most $n - 1$.  Let $NEAR(N,X)$ be the sum for $i = 1, \dots, n$ of the $\#NEAR(N,X,i)$, the number of integers in $NEAR(N,X,i)$.  
 
From $n > 1$ we have $0 < \#NEAR(N,X,i) < n$ which leads to $n \leq NEAR(N,X) \leq n(n -1)$ .  This upper bound is reached iff all the off-diagonal entries of $D(N,X)$ are the same nonnegative number.  When the entries $x_{ij}$ of $X$ satisfy no special relations (some would express this ``$X$ is chosen at random'') then $NEAR(N,X) = n$.  (Reason:  ``$x(j)$ and $x(j')$ are both $N$-nearest neighbors of $x(i)$'' is a special relation if $j$ and $j'$ are different.)  So, the set of $n \times k$ real matrices $X$ for which $NEAR(N,X) > n$ is a thin set in $nk$-dimensional real space, i.e. it contains no open $N$-ball of positive radius.  The symmetry of $D(N,X)$ implies $NEAR(N,X) \neq n(n-1) - 1$.

\vspace{0.5cm}

Question:  What values may $NEAR(N,X)$ assume?

\vspace{0.5cm}

\noindent ROBUSTNESS

\vspace{0.5cm}

Robustness will be defined in terms of nearest neighbors.  Throughout this section, we assume $n > 1$ to ensure that we have some nearest neighbors i.e. $NEAR(N,X) > 0$.

Let $X' \in$ $M(n \times (k + 1))$ have the same first $k$ columns as $X$ in $M(n \times k)$.  Let $x'(i)$ denote the $i$-th row of $X'$. 

\vspace{0.5cm}
            
Then $NEAR(N,X',i) = \{j \neq i : x'(j)$ is an $N$-nearest neighbor of $x'(i)\}$.

 \vspace{0.5cm}

\noindent In giving a meaning to``robustness'', (T23) says it is the proportion of nearest neighbors that remain unchanged in passing from $X$ to $X'$.  This can be interpreted as follows:  Let $m(i)$ be the number of $j$'s that are in both $NEAR(N,X,i)$ and $NEAR(N,X',i)$.  Then

$$ROB^{+}(N,X,X') = [m(1) + \cdots + m(n)]/NEAR(N,X)$$

\noindent is the proportion which (T23) calls the robustness of $X$.  However, we shall call this nonnegative rational quotient ``the $N$-robustness of $X$ with respect to $X'$'', thus emphasizing that it depends on $N$ and $X'$.  
It never exceeds $1$ because $m(i) \leq \#NEAR(N,X,i)$.  In fact, it is $< 1$ iff there is at least one $i$ such that $m(i) < \#NEAR(N,X,i)$.  If the final column of $X'$ is constant, e.g. all zeros, the distance matrices associated with $X$ and $X'$ are the same since $N$ satisfies (iv) and so the $N$-robustness of $X$ with respect to $X'$ will be $1$, the maximum possible value.  Ex6 will show that this robustness can be zero.  So, the only possible $N$-robustnesses of $X$ with respect to an $X'$, which is $X$ with one adjoined column, are the fractions

$$0/NEAR(N,X),\ 1/NEAR(N,X), \dots, NEAR(N,X)/NEAR(N,X),$$

\noindent and hence any limit of robustnesses along a sequence of members $X'$ of $M(n \times (k + 1))$ which are one-column extensions of $X$ is also in this set of fractions, and that sequence of fractions must ultimately become constant in order to converge.

\vspace{0.5cm}

\noindent Ex6:  Here is an example of $ROB^{+}(N_{p},X,X') = 0$, for all $p \in J$.
 
\noindent Let $X' = \begin{bmatrix} 2&50\\5&20\\1&10\end{bmatrix}$ and let $X$ and $Y$ denote respectively its first and second columns.  Since $N([1]) = 1$,

$$D(N,X) = \begin{bmatrix} 0&3&1\\3&0&4\\1&4&0\end{bmatrix}$$

\noindent and so the nearest neighbors correspond to positions $(1,3), (2,1), (3,1)$ in $D(N,X)$.  That is, these are the locations of the smallest off-diagonal entries in its rows.

$$D(N,Y) = \begin{bmatrix} 0&30&40\\30&0&10\\40&10&0\end{bmatrix}$$

Then

$$D(N_{\infty}, X') = \begin{bmatrix} 0&30&40\\30&0&10\\40&10&0\end{bmatrix},$$

\noindent each entry is the maximum of the entries in the same position in the two just preceding matrices.  So, the nearest neighbors correspond to positions $(1,2), (2,3), (3,2),$ which are disjoint from the positions above so $m(1) = m(2) = m(3) = 0$.  Hence $ROB^{+}(N_{\infty},X, X')$, the $N_{\infty}$-robustness of $X$ with respect to $X'$, is 0.

Let $p \in J$ be finite and we write down now not $D(N_{p},X')$ but $P$, the matrix whose entries are the $p$-th powers of the entries in $D(N_{p},X')$.

$$P =  \begin{bmatrix} 0 & 3^{p} + 30^{p} & 1^{p} + 40^{p} \\ 3^{p} + 30^{p} & 0 &4^{p} + 10^{p} \\ 1^{p} + 40^{p} & 4^{p} + 10^{p} & 0\end{bmatrix}.$$

If for all finite $p \in J$ the inequalities

$$4^{p} + 10^{p} < 3^{p} + 30^{p} < 1 +40^{p}$$

\noindent hold, then the $N_{p}$-nearest neighbors of $X'$ correspond to positions
$(1,2), (2,3), (3,2)$ of $D(N_{p},X')$ and as above $m(1) = m(2) = m(3) = 0$ and so $ROB^{+}(N_{p},X,X')= 0$ for all finite $p \in J$.

For the inequalities, observe:  $4^{p} + 10^{p} < 2(10^{p}) < 30^{p} < 3^{p} + 30^{p} < 2(30/40)^{p}40^{p} < 40^{p} < 1 + 40^{p}$.//

\vspace{0.5cm}

The definition of robustness we are using considers all possible augmentations of $X$ by one column.  A footnote in (T23) mentions another definition, found in (Temple 1982), which does not permit so much alteration of $X$:  Let $n_{j}$ be the number of rows of $X$, that is $i$'s in $\{1, ..., n\}$, such that $NEAR(N,X,i)$ changes when column $j$ is removed from $X$.  Here is the alternate definition of robustness of $X$ in $M(n \times k)$, which makes sense if $k > 1$:  
$$ROB^{-}(N,X) = 1 - (n_{1} + \cdots + n_{k})/nk.$$  

Here, only $k$ alterations of $X$ are allowed/taken into account, not infinitely many.  When $k > 1$ and all the rows of $X$ are the same then $D(N,Y)$ is $0(n \times n)$ whether $Y$ is $X$ or $X$ with a column removed and so $ROB^{-}(N,X) = 1 - 0/nk  =1$.

 \vspace{0.5cm}
 
And $ROB^{-}(N,X)$ too can vanish.

\vspace{0.5cm}

\noindent Ex7:  Let $p \in J$ be finite and $X$ be the first and $Y$ the second column of

$$Z = \begin{bmatrix} 1 & 0 \\ 0 & 0 \\ 0 & 1\end{bmatrix}.$$

\noindent Then

$$D(N_{p},Z) = \begin{bmatrix} 0 & 1 & 2^{1/p} \\ 1 & 0 & 1 \\ 2^{1/p} & 1 & 0 \end{bmatrix},$$

$$D(N_{p},Y) = \begin{bmatrix} 0 & 0 & 1 \\ 0 & 0 & 1 \\ 1 & 1 & 0 \end{bmatrix},$$

$$D(N_{p},X) = \begin{bmatrix} 0 & 1 & 1 \\ 1 & 0 & 0 \\ 1 & 0 & 0 \end{bmatrix}.$$

\noindent Then $n(1) = n(2) = 3$ and so the robustness is zero.  When $p = \infty$, then $2^{1/p} = 1$, so $ROB^{-}(N_{\infty}, Z) = 1/3$.

\vspace{0.5cm}

Exercise examples:  Let $Z$ denote the $3 \times 2$ matrix obtained by removing the zero first row of the matrix in Ex4 and let $X$ and $Y$ be the respective columns of the result.  Then $ROB^{-}(N_{p},Z) = 2/3$ for $p = 1, 2, \infty$. For note that the smallest off-diagonal entry in the first row of $D(N_{p},Z)$ is located in columns 2 and 3 for all three of the $p$ values.  The smallest in rows 2 and 3 are in columns 3 and 2 when $p = 1$ and in column 1 when $p = 2$ and $\infty$.

 \vspace{0.5cm}

\noindent Theorem:  Given a $p \in J$ and an $X \in M(n \times k)$ there is an $X' \in M(n \times (k + 1))$ which is $X$ with an additional last column such that

$$ROB^{+}(N_{p},X,X') \leq n/NEAR(N_{p},X).$$

Before our proof, which indicates how $X'$ can be constructed, let

$$T = \{(i,j) \in \{1, 2, ..., n\} \times \{1, 2, ..., n\} : i < j\}$$

\noindent and set $u(i,j) = 2j - 2i > 0$ for all $(i,j) \in T$.

\vspace{0.5cm}

\noindent Lemma:  If $(i,j)$ and $(r,s)$ are different elements of $T$, then $u(i,j)$ differs from $u(r,s)$.

\vspace{0.5cm}

\noindent Proof:  Were the two values of u to agree, there would be a largest integer $t$ such that both would be divisible by the $2^{t}$. Clearly $t = i$ and $r$.  So the lowest power terms are equal and hence the highest power ones are too, a contradiction. ///

\vspace{0.5cm}

Set $u(i,j) = u(j,i)$ when $j < i$ and $(j,i) \in T$.  Set $u(i,i) = 0$.  Then by the lemma the numbers $u(i,1), u(i,2), \dots , u(i,n)$ are all different.

\vspace{0.5cm}

\noindent Proof of the theorem:  Let $d(N_{p},i,j)$ (resp. $d'(N_{p},i,j)$) denote the $(i,j)$-th entry of $D(N_{p},X)$ (resp. $D(N_{p},X')$).

Let $y$ in $M(n \times 1)$ be the matrix whose $i$-th row $y(i) = 2i-1$.  Let $t > 0$ and let $ty$ be the last column of $X'$. Then

$$d'(N_{\infty},i,j) = \max\{d(N_{\infty},i,j), tu(i,j)\},$$

\noindent which is $tu(i,j)$ when $t$ is large enough.  Since the
$u(i,1), \dots , u(i,n)$ are all different there is a unique smallest positive one, say $u(i,j')$, and hence row $x'(j')$ is the unique nearest neighbor of $x'(i)$. We pick $t$ so big that this holds in every row of $D(N_{\infty},X')$.  Hence $NEAR(N_{\infty},X') = n$.

Now assume that $p \in J$ is finite and note the property that then

$$N_{p}(v) > N_{p}(w)$$

\noindent whenever $v,w \in M(1 \times k)$ and $v(1) \geq w(1), \dots , v(k) \geq w(k)$ and at least one of these inequalities is strict.

Again, the last column of $X'$ will be $ty$ but for a different $t > 0$.  Hence $d'(N_{p},i,j) > d(N_{p},i,j)$ for $j \in \{1, ..., i-1, i+1, ..., n\}$, by the norm inequality in the second sentence above.  While $x(i)$
has $\#NEAR(N_{p},X,i)$ nearest neighbors among the rows of $X$, row $x'(i)$ has only one nearest neighbor among the rows of $X'$, if we pick $t$ correctly:  Since the members of the $i$-th row of $[u(i,j)]$ are all different, and since $d(N_{p},i,j)$  for $j \in NEAR(N_{p},X,i)$ are all equal to the smallest off-diagonal enter in row $i$ of $D(N_{p},X)$, there is a unique $j \in NEAR(N_{p},X,i)$ such that $d'(N_{p},i,j)$ is the smallest positive number in

$$\{d'(N_{p},i,s) : s \in NEAR(N_{p},X,i)\},$$

\noindent and this $j$ is unaltered by changing $t > 0$.  To see that if we pick $t > 0$ small enough $x'(j)$ will be a nearest neighbor of $x'(i)$, suppose $m \notin NEAR(N_{p},X,i)$.  Then $d(N_{p},i,j) < d(N_{p},i,m)$ and as $t > 0$ approaches $0$, the number $d'(N_{p},i,j)$ approaches the left side of this inequality and  $d'(N_{p},i,m)$ approaches the right side.  Therefore, there is an $s > 0$ such that for all positive $t < s$ we have $d'(N_{p},i,j) < d'(N_{p},i,m)$.  We must do this for each $m \notin NEAR(N_{p},X,i)$ and for each we may be obliged to reduce the size of $s > 0$.  Now we know that $x'(j)$ is the only nearest neighbor of $x'(i)$ when $0 < t < s$.  Doing this for $i = 1, ..., n$ requires only a finite number of additional reductions in the size of $s > 0$, so then each row of $X'$ has just one nearest neighbor among the row of $X'$.  Thus $NEAR(N_{p},X') = n$.///

\vspace{0.5cm}

\noindent CONCORDANCE

\vspace{0.5cm}

Let $M, N$ be norms on $M(1)$ satisfying (iv) and (v) and let $X$ be in $M(n \times k)$.  Then the ``concordance of $D(M,X)$ with $D(N,X)$'' is the fraction

$$CORD(M,N,X) = \#\{i \in \{1, \dots, n\} : NEAR(M,X,i) = NEAR(N,X,i)\} /n.$$  

\noindent Clearly, $CORD(M,N,X) = CORD(N,M,X),$ but cf. Table 3 of (T23).
The range of $CORD(M,N,X)$ for $X \in M(n \times k)$ is within

$$\{0/n, 1/n, \dots, n/n\}.$$

Thus, by our theorem, so are any of its limits along a strictly increasing sequences of members of $M(n)^{1=} = M(n)$ and any converging sequence of its values is ultimately constant.
 
This definition does not match the one given at the end of the Introduction in (T23).  That definition is not as specific as ours, which is tied to a specified $X$.  Instead, (T23)'s definition seems to apply to each and every member of $M(1)$.    An additional problem is the mention of probability.  To mention that one must specify a sample space and a probability measure, or those which are implicitly intended must be arguably obvious.  The earliest drafts of (T23) included both definitions, perhaps thinking they were equivalent, though they aren't.  The ratio in our definition is a probability, each point of the sample space $\{1, \cdots, n\}$ having probability $1/n$.

\vspace{0.5cm}

\noindent EXAMPLES
\vspace{0.5cm}

Ex0:  When $n = 1$, we have $NEAR(M,X,1) = 0 = NEAR(N,X,1)$ so $CORD(M,N,X) = 1/1 = 1$.

Ex1:  When $n = 2$, $x(1)$ and $x(2)$ are $M$- and $N$-nearest neighbors of each other so $CORD(M,N,X) = 2/2 = 1$.

Ex2:  When $k = 1$, then $D(N,X)$ is the same for all $N$ satisfying (v) and so $NEAR(M,X,i) = NEAR(N,X,i)$ for all $i$, and so $CORD(M,N,X) = n/n = 1$.

Ex4: $X$ without its first row is  

$$Y = \begin{bmatrix}  2 &  0 \\  -1 & \sqrt{3} \\ -1 & -\sqrt{3} \end{bmatrix}.$$  

Then

$$D(N_{1},Y) = \begin{bmatrix} 0 & 3 + \sqrt{3} & 3 + \sqrt{3} \\ 3 + \sqrt{3} & 0 & 2\sqrt{3} \\ 3 + \sqrt{3} & 2\sqrt{3} & 0 \end{bmatrix}$$

$$D(N_{2}, Y) =  \begin{bmatrix} 0 & 2 \sqrt{3} & 2 \sqrt{3} \\ 2 \sqrt{3} & 0 & 2\sqrt{3} \\ 2 \sqrt{3} & 2\sqrt{3} & 0 \end{bmatrix}$$

\noindent Hence $CORD(N_{1},N_{2}, Y) = 1/3$.

Each of these examples shows that restricting attention to the $1$-norm and $2$-norm will not oblige a sequence of matrices in $M(n)$ with ever more columns, perhaps all the additional ones being $0(n \times 1)$, to drive the concordance to something near $0.570376001675023 = \delta = \exp(-\exp(-\gamma))$ where

$$\gamma = 0.5772156649015328606065$$

\noindent is the Euler-Mascheroni constant.  We are unaware of any opinions about whether $\delta$ is rational.  If it isn't and if $n$ is not allowed to vary as $k \to \infty$ then for no choices of $M$ and $N$ and $X$ in $M(n \times k)$ is $\delta$ the limit of $CORD(M,N,X)$ as $k \to \infty$ because each of these concordances is a rational number with denominator $n$.  Even if $\delta$ is rational it is not clear whether it is the limit of some $CORD(M,N,X)$s as the number of columns of $X$ increases but the number of rows stays fixed.  Working from the continued fraction expansion of $\delta$ shows that if $\delta = m/n$ for relatively prime positive integers $m$ and $n$ then the denominator $n$ must exceed $169229911$. (Thanks to L. N. Trefethen for getting MATLAB to compute this and the previous convergent.)  The previous convergent of the continued fraction expansion of $\delta$ has denominator $5382609$, so even if there is some roundoff error in the calculation, building an $X$ with more than $5$ million rows and the right number of nearest neighbors could be a challenge. If the computation of the expected distance to the nearest neighbor and the subsequent steps carried out in (T23) are justified and if they deliver the rational concordance defined above, that would provide a proof that $\delta$ is rational, a result/curiosity of some mathematical interest if it is not already known.

\vspace{0.5cm}

\noindent CONCERNS ABOUT (T23)'S TREATMENT OF CONCORDANCE

\vspace{0.5cm}
 
``We consider points forming a uniform random distribution in k-dimensional space with density $\lambda$ per unit volume, ...''.  Are these the $n$ points mentioned in (T23) in the first sentence of the Abstract and of the third paragraph of the Introduction, that is, the $n$ rows of $X$?  If so, they are all contained in some origin-centered ball and outside it their density seems to be zero.  How was their/any density supposed to be measured?  That is not explained.  We take ``uniform'' to mean that each point has the same probability measure, namely $1/n$, not that the points are in some sense evenly spread around/distributed.  However, the concept of ``evenly distributed'' does show up in (T23) in the paragraph before its equation (3).

Were the set of ``points'' infinite, this $1/n$ would not work and in one-dimensional space $-1$ would not have a nearest neighbor in the set $\{1, 1/2, 1/3, ...\}$.  What does randomly distributed mean?  Not arranged in any special pattern?  Can one define density and find an example of such a finite set meeting all these requirements even in $1$-dimensional space?  Is this challenge meetable if the set is infinite?  The cited informal discussion/result in Feller's book (p. 158) does mention (but not define) density, so how does (T23)'s $\lambda$ relate to Feller's?  Perhaps we should just think of $\lambda$ as a positive constant and stop trying to interpret it as density of the set of points.

Feller's result (for which he does not include a formal proof) seems to have its dependence on the number of points subsumed in $\lambda$.  At one point in (T23) $\lambda$ is set equal to $1$.  How is that justified? Doesn't that imply that the set is topologically dense? (If not, it omits an open interval, and so can't have density $1$ on any interval containing the omitted one.)  But then, some limit points could lack any nearest neighbor, e.g. the limit point $0$ of the harmonic sequence.

Relying on (Feller p. 158) (T23) finds that the average/expected distance from $0$ to the closest member of all the random uniformly distributed subsets of $k$-dimensional space with density $\lambda$ is finite.  For comparison we now ask just in $1$-dimensional space what the average/expected distance from $0$ to the closest member of $n$ points in the space is.  We will view all these sets of n points as equally probable.  (We even allow some of the points to equal each other, as that costs us nothing.)  Let $L > 0$ and use the uniform = constant distribution on the sample space which is the real interval $[-L,L]$.  We begin with singleton sets and uniform distribution $f(x) = 1/(2L)$ for $x \in [-L,L]$.  The expected value is the integral over the interval sample space of $(2L)^{-1}|x|dx, |x|$ being the distance of $x$ from $0$.  The positive and negative halves of the sample space contribute the same amount, so the result is $2(2L)^{-1}[(L^{2}/2) - (0^{2}/2)] = L/2.$  Thus the expectation increases without bound as $L$ increases without bound.  This is unsurprising because the uniform distribution treats all the $x$ values the same.  As $|x|$ wanders more widely its average position increases.  (T23) obtains a finite value over the whole real line because the distribution used there goes to zero at infinity (and $0$) prejudicing the calculation in favor of the smaller values of $|x|$.

We leave as an exercise the case in which the sets are not singletons, but contain two points.  

Now suppose the sets are triples still in $1$-dimensional space.  Here the uniform distribution is $f(x,y,z) = (2L)^{-3}$ for $(x, y, z)$ in the origin centered $L_{\infty}$-ball of radius $L$, that is the origin-centered cube with edge length $2L$ and eight vertices $(\pm L,\pm L, \pm L)$.  Consider the points $(x,y,z)$ of the sample space which are also in the first octant and have $x$ a lower bound of $y$ and of $z$.  Their contribution to the expected value is

$$E = (2L)^{-3}\int_{0}^{L}\int_{x}^{L}\int_{x}^{L} x dzdydx.$$

\noindent There is an equal contribution to the expectation by the first octant points of the sample space having $y$ (resp. $z$) a lower bound of $x$ and of $z$ (resp. $y$).  And since there are eight octants, the expectation is $24E = L/4$.

The easy guess is that the expectation computed for sets with $n$ points will be $L/(n+1)$ and will approach infinity as $L$ does.  The finite expectation in (T23) arises from relying on (Feller 1968) and favoring some points over others.  Our approach treats all points the same.  The two methods disagree even on a one- dimensional space.

\vspace{0.5cm}

\noindent THE CORRELATION COEFFICIENT OF TWO DISTANCE MATRICES

\vspace{0.5cm}

In order to define the correlation between two random variables, we need a sample space $\Omega$ on which they are defined and a probability measure defined on $\Omega$.  One natural choice when working on $n \times n$ distance matrices is

$$\Omega = \{1, \dots, n\} \times \{1, \dots, n\}$$

\noindent with $1/n^{2}$ the measure of each of its points.  Since distance matrices are symmetric and zero along their diagonals, another natural choice is
 
$$\Omega = T = \{ (i,j) \in \{1, \dots, n\} \times \{1, \dots, n\} : i < j\}$$

\noindent with $2/n(n - 1)$ the measure of each of its points.  Let $D$ be a distance matrix and $S$ denote the sum of the entries of $D$ which are above its main diagonal.  Then the expectation $\mathbb{E}(D)$ of $D$ is $2S/n^{2} = (1/n)(2S/n)$ using the first sample space and the easily comparable $2S/n(n-1) = (1/(n - 1))(2S/n)$ using the second.  We use the first in sequel.

Recall that if $A = [a_{ij}]$ and $B = [b_{ij}]$ are matrices of the same size then their Hadamard product $A \circ B = [a_{ij}b_{ij}]$ also has the same size.  Let $M, N$ be two norms satisfying (iv) and (v), and $X$ be in $M(n \times k)$, as usual.  Then the correlation coefficient of $D(M,X)$ and $D(N,X)$ will be denoted $\rho(M,N,X)$ and (with var = variance and cov = covariance) is

$$\frac{cov(D(M,X),D(N,X))}{\sqrt{var(D(M,X))\ var(D(N,X))}},$$

\noindent which is known to equal

$$\frac{\mathbb{E}(D(M,X) \circ D(N,X)) - \mathbb{E}(D(M,X)\mathbb{E}(D(N,X)}{\sqrt{(\mathbb{E}(D(M,X) \circ D(M,X)) - (\mathbb{E}[D(M,X)])^{2})(\mathbb{E}(D(N,X) \circ D(N,X)) - (\mathbb{E}[D(N,X)])^{2})}}.$$

\noindent Since $A \circ B = B \circ A$ it follows that $\rho(M,N,X) = \rho(N,M,X)$.  This theoretical symmetry is not perfectly reflected in Table 3 in (T23).  More information about how the numbers in the table were derived might explain this.

\vspace{0.5cm}

\noindent EXAMPLES

\vspace{0.5cm}

Ex0:  If $n = 1$ then $var(D(M,X)) = var(D(N,X)) = var([0]) = 0$ so $\rho(M,N,X)$ is not defined.

Ex1:  If $n = 2$ and the rows of $X$ are the same then $D(M,X) = 0(2 \times 2) = D(N,X)$, which have variance zero, so $\rho(M,N,X)$ is not defined.  If the rows are different then $w = x(2) - x(1) \neq 0$ and $a = M(w) > 0$ and $b = N(w) > 0$.  So $(b/a)D(M,X) = (b/a)\begin{bmatrix} 0 & a \\ a & 0 \end{bmatrix} = D(N,X)$.  Hence $\rho(M,N,X) = 1$ because the two distance matrices are scalar multiples of each other.

Ex2:  $X = [x_{11} \  \dots \  x_{n1}]^{T} \neq 0(n x 1)$ and so $D(N,X) =  [|x_{j1} - x_{i1}|]$ for all norms $N$ satisfying (v).  So $\rho(M,N,X) = 1$, if norm $M$ satisfies (v) too.

Ex3':  If in Ex3 above $0(1 \times k) \neq w$ for any integer $k > 0$, then $D(M,X)$ is a nonzero scalar multiple of $D(N,X)$, so $\rho(M,N,X) = 1$.

Ex8:  Let $X$ be the first column of

$$Y = \begin{bmatrix} 1&1\\ 2&0\\ 3&0 \end{bmatrix}.$$

\noindent In order to avoid the computationally clumsy square root in
$N_{2}(v)$, (T23) considers in its place $L(v) = (N_{2}(v))^{2}$, which, since it violates (ii) and (iii) is not a norm, but satisfies (iv) and (v).  For all $p \in J$

$$D(N_{p},X) = \begin{bmatrix} 0 & 1 & 2 \\ 1 & 0 & 1 \\ 2 & 1 & 0 \end{bmatrix},$$

$$D(L, X) = D(N_{p},X) \circ D(N_{p},X) =  \begin{bmatrix} 0 & 1 & 4 \\ 1 & 0 & 1 \\ 4 & 1 & 0 \end{bmatrix},$$

$\mathbb{E}(D(N_{p},X)) = (2/9)(4),$

$\mathbb{E}(D(N_{p},X) \circ D(N_{p},X)) = (2/9)(6),$

$\mathbb{E}(D(L,X)) = (2/9)(6),$

$\mathbb{E}(D(L,X) \circ D(L,X)) = (2/9)(18),$

$\mathbb{E}(D(N_{p},X) \circ D(L,X)) = (2/9)(10),$

\noindent then $\rho(N_{p},L,X) = 7/\sqrt{55} = 0.94388$.

Now we turn to $Y$ and restrict $p$ to $1$ and $\infty$.  

$$D(N_{1},Y) = \begin{bmatrix} 0 & 2 & 3\\ 2 & 0 & 1 \\ 3 & 1 & 0 \end{bmatrix},
$$

$$D(L,Y) = \begin{bmatrix} 0 & 2 & 5\\ 2 & 0 & 1 \\ 5 & 1 & 0 \end{bmatrix},$$

$$D(N_{\infty},Y) = \begin{bmatrix} 0 & 1 & 2\\ 1 & 0 & 1 \\ 2 & 1 & 0 \end{bmatrix}.$$

$\mathbb{E}(D(N_{1},Y)) = (2/9)(6),$

$\mathbb{E}(D(N_{1},Y) \circ D(N_{1},Y)) = (2/9)(14),$

$\mathbb{E}(D(L,Y)) = (2/9)(8),$

$\mathbb{E}(D(L,Y) \circ D(L,Y)) = (2/9)(30),$

$\mathbb{E}(D(N_{1},Y) \circ D(L,Y)) = (2/9)(20),$

\noindent then $\rho(N_{1},L,Y) = 14/\sqrt{213} = 0.959264.$

\vspace{0.5cm}

$\mathbb{E}(D(N_{\infty},Y)) = (2/9)(4),$

$\mathbb{E}(D(N_{\infty},Y) \circ D(N_{\infty},Y)) = (2/9)(6),$

$\mathbb{E}(D(N_{\infty},Y) \circ D(L,Y)) = (2/9)(13),$

then $\rho(N_{\infty},L,Y) = 53/(2\sqrt{781}) = 0.948245$.

\vspace{0.5cm}

\noindent REMARKS
\vspace{0.5cm}

The claim that $\rho(N_{1},L,W)$ for $W$ in $M(n \times k)$ is independent of $k$ (which is true, as explained far above, if for example we increase the number of columns of $W$ by adjoining only constant columns to it) can be interpreted two ways:  one that it applies to matrices built up from $W$ by adjoining additional columns.  This will not work if adjoining just any column is allowed as one sees by checking that $\rho(N_{1},L,Y) = 14/\sqrt{213} = 0.959264 \neq \rho(N_{p},L,X) = 0.94388$.  And this fact also faults the other interpretation of "independent of $k$", namely that

$$\rho(N_{1},L,Z) = \rho(N_{1},L,W)$$

\noindent for every $Z$ with $n$ rows and more columns than $W$.

Of course, by giving $Y$ additional columns all of which are constant, the value of $\rho(L,N_{\infty},Y) = \rho(N_{\infty},L,Y)$ will not change, in particular not decrease as the number of columns increases and not approach $0$.  Cf. Table 3 of (T23).

Ex9:  Let $r = \sqrt{3}$ and set

$$X = \begin{bmatrix} 2 & 0 \\ -1 & r \\ -1 & -r \end{bmatrix}$$

\noindent cf. Ex4 without its first row.  We have

$$D(N_{1},X) = \begin{bmatrix} 0 & 3 + r & 3 + r\\ 3 + r & 0 & 2r \\ 3 + r & 2r & 0 \end{bmatrix}.$$

\noindent $\mathbb{E}(D(N_{1},X))= (2/9)(6 +4r)$ and $var(D(N_{1},X)) = 2/9(36 + 12r) - (2/9(2/9)(6 + 4r)(6 + 4r) = (8/27)(13 + r).$

$$D(N_{2},X) = \begin{bmatrix} 0 & 2r & 2r\\ 2r & 0 & 2r \\ 2r & 2r & 0 \end{bmatrix}.$$

\noindent $\mathbb{E}(D(N_{2},X) = (2/9)(6r)$ and\\ $var(D(N_{2},X)) = (2/9)(36) -(2/9)(6r)(2/9)(6r) =(4/81)(162 - (6r)(6r)) = (4/81)(54) = 8/3.$

$$D(N_{\infty},X) = \begin{bmatrix} 0 & 3 & 3 \\ 3 & 0 & 2r \\ 3 & 2r & 0 \end{bmatrix}.$$

\noindent $\mathbb{E}(D(N_{\infty},X) = (4/9)(3 + r)$ and $var(D(N_{\infty},X)) = 2/9(30) - (2/9(2/9)(6 + 2r)(6 + 2r) = (4/27)(29 - 8r)$

\noindent $cov(D(N_{1},X),D(N_{2},X)) = (8/9)(2 + r)$

\noindent $cov(D(N_{1},X),D(N_{\infty},X)) = (4/27)(25 - 3r)$

\noindent $cov(D(N_{2},X),D(N_{\infty},X)) = (8/9)(1 + r)$

$\rho(N_{1},N_{2},X) = 0.972335$

$\rho(N_{1},N_{\infty},X) = 0.9375373$

$\rho(N_{2},N_{\infty},X) = 0.9928629$

\vspace{0.5cm}

\noindent REMARKS
\begin{enumerate}
\item As we remarked far above, these numbers (which are not guaranteed accurate to all the decimal places displayed) will not change if $X$ is augmented with any number of constant columns.  So when $p = 1 ,2$ we know that $\rho(N_{p},N_{\infty},X)$ need not approach 0 as the number $k$ of columns of $X$ increases.  Cf. Table 3 of (T23).

\item Every correlation we have computed is larger than $9/10$.  The maximum possible is $1$.  What is the minimum possible correlation?  That is, given norms $M$ and $N$ find the greatest lower bound of $\rho(M,N,X)$ as $X$ ranges over $M(n \times k)$.  Is there an $X$ at which the correlation assumes this lower bound?  I.e. Is it a minimum?  Our examples above have settled this question when $n = 1$, when $n = 2$, and when $k = 1$, so one might start with $X$s in $M(3 \times 2)$ and $M$ and $N$ different $p$-norms with $p \in \{1, 2, \infty\}$.

\item Does $ROB^{+}(N,X,X')$ increase with the number $k$ of columns in $X$
while the concordance $CORD(N_{1},N_{2},X)$ decreases with $k$ as illustrated in Fig. 3 of (T23)?  Not necessarily.  Our trick of augmenting $X$ with constant columns leaves both these numbers fixed.
\end{enumerate}

\vspace{0.5cm}
\noindent APPENDIX

\vspace{0.5cm}

If a set $B$ in $k$-dimensional real space contains an open origin-centered ball of positive radius and has volume $V_{0}$, and then each of the points in it is multiplied by a scaling factor $r > 0$ the set which results has volume $V_{0}r^{k}$.  To see this is reasonable, just imagine that $B$ is an assemblage of cubes, some tiny, and observe that scaling a cube changes its volume this way and hence so does scaling the assemblage.

If
$$f(r) = d/dr(\exp[-\lambda V_{0}r^{k}]) = -k\lambda V_{0}r^{k - 1}\exp[-\lambda V_{0}r^{k}]$$

\noindent is taken to be the probability density function of $r$, as it is in (T23), then

$$\mathbb{E}(r) =  \int_{0}^{\infty}rf(r)dr = -\lambda V_{0}\int_{0}^{\infty}kr^{k}\exp[-\lambda V_{0}r^{k}]dr = \Gamma(1 + 1/k)/(\lambda V_{0})^{1/k}.$$

\noindent (T23) cites one place this integral can be found.  It can also be viewed as the Laplace transform of $u^{1/k} = r$ evaluated at $\lambda V_{0}$.

With this value of $r$, the volume is

$$V_{0}(\Gamma(1 + 1/k)^{k}/(\lambda V_{0}) = [\Gamma(1 + 1/k)]^{k}/\lambda.$$

\noindent This volume is then independent of the original set $B$ and its volume $V_{0},$ which may surprise, and we see that the calculation does not need to be repeated, as (T23) does, when the shape or volume of $B$ changes.

\vspace{0.5cm}

\noindent {\bf Funding:}  The author was not funded by anyone.\\
{\bf Conflict of Interest:}  The author was not associated with any conflicting interests.\\
{\bf Ethical Conduct:}  The author understands that he is fully ethically compliant.\\
{\bf Data Availability:}  No data were involved in this work.\\
{\bf Acknowledgement:}  Professor Mark Mills shepherded my LaTeX file through this journal's
electronically demanding submission process.  I could not have done it myself.  Thank you,
Professor Mills.

\vspace{0.5cm}

\noindent REFERENCES

\vspace{0.5cm}
\begin{enumerate}
\item FELLER, W. (1968). Introduction to probability theory and its applications (3rd edn), vol. 1.  Wiley.  xviii +509 pp.

\item MAZUR, S. \& ULAM, S. (1932). Sur les transformations isom$\acute{e}$triques d'espaces vectoriels norm$\acute{e}$s.  Comptes Rendus de l'Acad. des Sc., 194, 946-948.

\item TEMPLE, J. T. (1982). An empirical study of robustness of nearest-neighbor relations in numerical taxonomy. Mathematical Geology, 14, 675-678.

\item TEMPLE, J. T. (2023). Characteristics of distance matrices based on Euclidean, Manhattan and Hausdorff coefficients.  Journal of Classification, 40, 214-232.
\end{enumerate}

\end{document}